\documentclass[10pt]{article}

\usepackage{amsmath}         
\usepackage{amssymb} 
\usepackage{amsthm}  
\usepackage{enumerate}
\usepackage{stmaryrd}
\usepackage{url} 
\usepackage{mathrsfs}       
\usepackage[all]{xy}
\usepackage{graphicx}

\usepackage{times}

\theoremstyle{plain}

\theoremstyle{definition}

\def\Ro{R_{\rm o}}

\def\M{{\cal M}}

\def\G{\varGamma}

\def\nx{\raise 1pt\hbox{$\bigcirc$}} 

\def\T{{{\mit\Theta}}}

\def\rtb{R^\T_{ \scriptscriptstyle\Box}}

\def\\{\vert}

\title{A simplified temporal logic completeness proof }
\author{Robert Goldblatt\thanks{School of Mathematics and Statistics,
Victoria University, Wellington, New Zealand.}}

\date{}

\begin{document}

\maketitle

\begin{abstract}

A temporal logic discussed in the author's book \emph{Logics of Time and Computation} is given a  simpler proof of semantic completeness, eliminating some flaws in the original.
\end{abstract}

\bigskip
This note assumes familiarity with the content of Chapter 9 of [1], which we   refer to frequently by citing definitions and results and giving page references.

Chapter 9 gives a  semantic completeness proof  for a propositional temporal logic  $\T$ whose formulas are evaluated on infinite time sequences $$\sigma_0,\sigma_1,\dots\sigma_n,\dots\dots$$
where each $\sigma_n$ is thought of as  the \emph{state} produced at time $n$ by the execution of some computation. 

The completeness argument on pages 91--98   makes use of the relation $\Ro^c$ defined on page 95, and the schema $Dum$. This builds on the use of $R^c$ and $Dum$ in the completeness theorem in Chapter 8 for the  $\Box$-only logic $S4.3Dum$. 
But in the second paragraph of page 96 there is an unjustified assertion about a finite sequence of simple $\Ro^c$-clusters giving rise to an $\Ro$-list. On the same page, the assertion labelled (2), which refers to  $\Ro^c$-clusters, also lacks justification.

However those assertions are not needed. Moreover $\Ro^c$ and $Dum$ themselves are not needed.
The formal language studied in Chapter 9 is richer than that of Chapter 8. It includes the $\nx$ (``next'') modality, which appears in axioms and theorems containing the other two non-Boolean connectives $\Box$ (``henceforth") and $\mathcal{U}$ (``until"). As well as allowing the $\nx$-free schema $Dum$ to become derivable [1, Theorem 9.2(6)], they facilitate a simpler approach to  proving completeness  for $\T$ (see the use of $B\land\nx\Box B\to \Box B$ in the proof of the Lemma below), an approach that does not involve $\Ro^c$ or $Dum$ at all.

To explain this,  recall that the completeness proof begins with a non-$\T$-theorem $A$ and a $\T$-maximal (i.e.\ maximally $\T$-consistent) set $s_A$ of formulae with $A\notin s_A$. The goal is to show that $A$ is falsifiable in a model on a state sequence. 

Let $S^\T$ be the set of $\T$-maximal sets of formulae. This carries relations $\rtb$ and $\Ro^\T$, where $s\rtb t$ iff
$\{B:\Box B\in s\}\subseteq t$, and $\Ro^\T$ iff  $\{B:\nx B\in s\}\subseteq t$.
Put $S=\{u\in S^\T: s_A \rtb u\}$.
$\rtb$ is reflexive and transitive, and its restriction to $S$ is connected  [1, page 91].

Now let $\G$ be the finite set of formulas defined from the subformulas of $A$ as on [1, page 92],
and consider the  structure  $(S_\G,\Ro,\Ro^*)$
introduced on [1, page 95]. Here $S_\G$ is the finite set obtained by $\G$-filtration of $S$, so that
$S_\G=\{\\s\\:s\in S\}$, where $\\s\\=\{t\in S: s\cap \G=t\cap\G\}$.
$\Ro$ is the smallest $\G$-filtration of $\Ro^\T$ [1, page 94], and $\Ro^*$ is the ancestral (reflexive transitive closure) of $\Ro$.  By the Ancestral Lemma 9.8, $\Ro^*$ is a $\G$-filtration of $\rtb$  \,[1, page 93].

An \emph{$\Ro$-list from $x$ to $y$} is a finite sequence $x=x_0,\dots,x_n=y$ such that $x_i\Ro x_{i+1}$ for all $i<n$.
In general we have $x\Ro^*y$ iff there is an $\Ro$-list from $x$ to $y$ of some finite length $n\geq 0$ \,[1, pages 9--10]. 

As noted on [1, page 95], $(S_\G,\Ro,\Ro^*)$ consists of a finite sequence of  $\Ro^*$-clusters. These are the equivalence classes of the  relation  $\{(x,y): x\Ro^*y\Ro^*x\}$ [1, page 66].
They are the only kind of cluster we consider. If $x$ and $y$ belong to the same $\Ro^*$-cluster then any $\Ro$-list from $x$ to $y$ is contained within that cluster.
We can use this to  ``unwind'' each $\Ro^*$-cluster  into a finite $\Ro$-list. By concatenating  lists we can create one that starts at any given point $x$, includes all members of the $\Ro^*$-cluster of $x$, and ends at any prescribed point in that cluster, including $x$ itself. We can then extend the list to any point in any later $\Ro^*$-cluster.

So we proceed ``from left to right'', unwinding each  $\Ro^*$-cluster in succession  until we have a single such $\Ro$-list that leads from the first $\Ro^*$-cluster, includes all of the points in the non-last $\Ro^*$-clusters, and ends at some point $x$ in the last $\Ro^*$-cluster. Then we unwind the last cluster into an $\Ro$-list that starts and finishes at $x$ and includes all members of the last cluster, and repeat that last list ad infinitum.

The result is an infinite sequence $\sigma:\omega\to S_\G$ that includes all points of $S_\G$. Its essential properties are that $\sigma_j\Ro\sigma_{j+1}$ in general, hence $\sigma_j\Ro^*\sigma_k$ whenever $k\geq j$, and that each member of the last $\Ro^*$-cluster is repeated infinitely often in $\sigma$.

Now for an atomic formula $p\in\G$, we can well-define $V_\G(p)=\{\\s\\ \in S_\G: p\in s\}$.  This allows us to form the model                  
${\M}=(S_\G,\sigma,V_\G)$ 
for the language of formulae whose atomic subformulae come from $\G$.
For this we wish to prove 

\medskip
{\bf Theorem 9.13}:   if $B\in\G$, then for any $j\in \omega$ and $s\in\sigma_j$,
\begin{equation}
B\in s \quad{\rm iff}\quad \M \models_j B.
\end{equation}
A fundamental fact that we use is  $Fun$-Lemma 9.9 [1, page 94], proven using the axiom
 $Fun$:   $\nx\neg A \leftrightarrow \neg\nx A$.
It states that
if $\nx B \in \G$ and $s \in S$, then
$$
\text{
 $\nx B \in s$  \enspace iff \enspace
$\forall t(\\s\\\Ro\\t\\$ implies $B \in t$) \enspace iff \enspace
$\exists t(\\s\\\Ro\\t\\$ and $B \in t$).}
$$
The proof of (1)  is by induction on the formation of $B$. All cases hold just as in the proof of Theorem 9.13 on pages 96-98 of [1], except for the inductive case of a formula $\Box B\in\G$ having $\Box B\notin s$  (third paragraph of page 97). For this case we now provide the following Lemma, which states that there is some $k\geq j$ with $B\notin t$ for any $t\in\sigma_k$. Hence $ \M\not \models_k B$  by induction hypothesis on $B$, implying 
$ \M\not \models_j \Box B$ as required in this case.

\medskip\noindent
{\bf Lemma.} {\it
Let\/ $\Box B\in \G$ and $j\in\omega$. If\/ $\Box B\notin s$, where $s\in\sigma_j$, then there exists some
$k\geq j$ such that $B\notin t$ for all $t\in\sigma_k$.}

\medskip\noindent
{\it Proof.}
Suppose, for the sake of contradiction, that  $\Box B\notin s\in\sigma_j$, but there is no such $k$.
Then for all $k\geq j$, there exists some $t\in\sigma_k$ such that $B\in t$.
Now as  $B\in\G$, $B$ belongs to one member of  $\sigma_k$ iff it belongs to all members of $\sigma_k$. Hence

\begin{equation}
\text{ for all $k\geq j$, $B\in  t$ for all $t\in\sigma_k$.}
\end{equation}

\noindent
In particular, $B\in s$ as $s\in\sigma_j$.   Now $\Box B\notin s$ by hypothesis, and $B\land\nx\Box B\to \Box B$ is $\varTheta$-derivable (Theorem 9.2(4)), so we must have  $\nx\Box B \notin s$. But  $\nx\Box B \in\G$ and $|s|=\sigma_j\Ro\sigma_{j+1}$, so 
$\Box B\notin  t$ for all $t\in\sigma_{j+1}$ by  $Fun$-Lemma 9.9.
 
 We can then repeat this argument with $j+1$ in place of $j$: for all $t\in\sigma_{j+1}$ we have $\Box B\notin  t$ as just shown, while $B\in t$ by (2),  so $\nx\Box B \notin t$, hence $\Box B\notin  t'$ for all $t'\in\sigma_{j+2}$. Making this argument into an  induction on $k$,  we obtain  that 
 \begin{equation}
\text{ for all $k\geq j$, $\Box B\notin  t$ for all $t\in\sigma_k$.}
\end{equation}

\noindent
Now choose some $k\geq j$ such that $\sigma_k$ is in the last  $\Ro^*$-cluster, and take any
$t\in\sigma_k$. Then $\Box B\notin t$ by (3), so in the canonical model we have some $u$ with 
$t\rtb u$ and $B\notin u$.  Since $\Ro^*$ is a $\G$-filtration of $\rtb$ [1, page 93], from $t\rtb u$ we get $\sigma_k=|t| \Ro^*|u|$, so $|u|$ is also in the last $\Ro^*$-cluster.  Since the members of that cluster are repeated infinitely often in the sequence $\sigma$, we have $|u|=\sigma_l$ for some $l\geq k$. Then $l\geq j$ and $B\notin u\in\sigma_l$.
However, this  contradicts (2), and so the Lemma is proved, yielding Theorem 9.13 for this construction as indicated.

\medskip
The completeness proof for $\varTheta$ then follows as on [1, page 98]: since the sequence $\sigma$ includes all points of $S_\G$, there is a $j\in\omega$ with $\\s_A\\= \sigma_j$. As $A\in\G$ and $A\notin s_A\in\sigma_j$, Theorem 9.13 then gives  $\M\not\models_j A$, fulfilling our goal of showing that $A$ is falsifiable in a model on a state sequence.

\bigskip
\subsection*{Reference}

\begin{description}
\item{[1]}
Robert Goldblatt.
{\em Logics of Time and Computation}.
CSLI Lecture Notes No. 7. CSLI Publications, Stanford University,
  second edition, 1992.
\end{description}

\end{document}